\documentclass[12pt]{amsart}
\usepackage{amssymb}
\usepackage{enumitem}
\usepackage{graphicx}

\makeatletter
\@namedef{subjclassname@2020}{%
  \textup{2020} Mathematics Subject Classification}
\makeatother

\usepackage[T1]{fontenc}
\newtheorem{theorem}{Theorem}[section]

\newtheorem{lemma}[theorem]{Lemma}

\newtheorem{problem}[theorem]{Problem}

\theoremstyle{definition}

\numberwithin{equation}{section}
\def\pmod #1{\ ({\rm{mod}}\ #1)}

\def\bg{\bigg}
\def\({\bg(}
\def\){\bg)}

\def\Arg{{\rm Arg}}

\def\ve{\varepsilon}

\theoremstyle{plain}

\begin{document}

\baselineskip=17pt
\hbox{accepted version of 201108-Wang for Colloq. Math.}
\medskip

\title[Products of quadratic residues and related identities]{Products of quadratic residues and
related identities}
\author[Hai-Liang Wu and Li-Yuan Wang] {Hai-Liang Wu and Li-Yuan Wang}

\address {(Hai-Liang Wu)  School of Science, Nanjing University of Posts and Telecommunications,
Nanjing 210023, People's Republic of China}
\email{whl.math@smail.nju.edu.cn}

\address
{(Li-Yuan Wang, corresponding author) School of Physical and Mathematical Sciences, Nanjing Tech
University, Nanjing 211816, People's Republic of China}
\email {wly@smail.nju.edu.cn}

\date{}

\begin{abstract}
In this paper we study products of quadratic residues modulo odd primes and  prove some identities
involving quadratic residues.
For instance, let $p$ be an odd prime. We prove that  if $p\equiv5\pmod8$, then
$$\prod_{0<x<p/2,(\frac{x}{p})=1}x\equiv(-1)^{1+r}\pmod p,$$
where $(\frac{\cdot}{p})$ is the Legendre symbol and $r$ is the number of $4$-th power residues
modulo $p$ in the interval $(0,p/2)$.
Our work  involves class number formula, quartic Gauss sums, Stickelberger's congruence and values
of Dirichlet L-series at negative integers.
\end{abstract}

\subjclass[2020]{Primary 11A15; Secondary 11R11, 11R18.}

\keywords{Quadratic residues, cyclotomic fields, Gauss sums.}

\maketitle
\section{Introduction}

 The study of quadratic residues modulo odd primes is one of the most classical topics in number
 theory, and it has deep relations with other areas in number theory such as Gauss sums and
 permutations over finite fields. For example, Gauss in 1811 determined the explicit value of
 quadratic Gauss sums
(cf. \cite[Chapter 5]{IR}), i.e.,
 \begin{equation*}
   \tau_p:=\sum_{k=1}^{p-1}\(\frac{k}{p}\)e^{2\pi k{\bf i}/p}=\sqrt{(-1)^{(p-1)/2}p},
 \end{equation*}
 where $(\frac{\cdot}{p})$ denotes the Legendre symbol and ${\bf i}$ is the $4$-th primitive root of
 unity with $\Arg({\bf i})\equiv \pi/2\pmod{2\pi\mathbb{Z}}$ (where $\Arg(z)$ denotes the argument of
 a complex number $z$). We know that Gauss sums have many applications in number theory and in the
 study of finite fields. As another application of quadratic residues, Sun \cite{Sun ffa} investigated
 many permutations over finite fields involving squares in finite fields and in the same year he
 \cite{Sun ffa1} studied some determinants concerning the Legendre symbol.

 Let $p$ be a prime with $p \equiv 1 \pmod 4,$ and let $\zeta_p=e^{2 \pi{\bf i}/p}$ be the primitive
 $p$-th root of unity.   Sun \cite{Sun ffa} proved that
 $$
\prod_{k=1}^{(p-1) / 2}\left(1-\zeta_p^{ k^{2}}\right)=\begin{cases}\sqrt{p} \varepsilon_{p}^{-
h(p)}&\mbox{if}\
	p\equiv 1\pmod 4,\\(-1)^{(h(-p)+1) / 2} {\bf i}\sqrt{p} &\mbox{if}\ p\equiv 3\pmod4,\end{cases}
$$
where $\varepsilon_{p}>1$ and $h(p)$ are the fundamental unit and the class number of the real
quadratic field $\mathbb{Q}(\sqrt{p})$ respectively, and $h(-p)$ is the class number of  the imaginary
quadratic field $\mathbb{Q}(\sqrt{-p})$.
In the same paper he also proved the following results:
\begin{itemize}
\item If $p \equiv 1(\bmod 4),$ then
$$
\prod_{0<j<k<p/2}\left(\zeta_p^{ j^{2}}-\zeta_p^{k^{2}}\right)^{2}=(-1)^{(p-1) / 4} p^{(p-3) / 4}
\varepsilon_{p}^{  h(p)},
$$
\item If $p \equiv 3(\bmod 4),$  then
$$\prod_{0<j<k<p/2}\left(\zeta_p^{j^{2}}-\zeta_p^{ k^{2}}\right)=\begin{cases}(-p)^{(p-3) /
8}&\mbox{if}\ 8\mid p-3,\\(-1)^{\frac{p+1}{8}+\frac{h(-p)-1}{2}}  p^{(p-3) / 8} &\mbox{if}\ 8\mid
p-7.\end{cases} $$
\end{itemize}
Let $\# S$ denote the cardinality of a set $S$. Later Petrov and Sun \cite{Sun era} further obtained
the following results:
\begin{itemize}

\item  If $p \equiv 1(\bmod 8),$ then
$$
\prod_{0< j<k <p/2}\left(\zeta_p^{ j^{2}}+\zeta_p^{ k^{2}}\right)=(-1)^{\#\left\{1 \leqslant
k<\frac{p}{4}:\left(\frac{k}{p}\right)=-1\right\}},
$$
\item  If $p \equiv 5$ (mod 8 ),  then
$$(-1)^{\#\left\{1 \leqslant k<\frac{p}{4}:\left(\frac{k}{p}\right)=-1\right\}} \prod_{0<j<k
<p/2}\left(\zeta_p^{ j^{2}}+\zeta_p^{ k^{2}}\right)=  \varepsilon_{p}^{-  h(p)}.$$
\end{itemize}

Now we present some earlier results on sums of quadratic residues. For any prime $p>3$, let
$$\mathcal{R}:=\bigg\{x\in\mathbb{Z}: \(\frac{x}{p}\)=1\bigg\}\ \text{and}\
\mathcal{N}:=\bigg\{x\in\mathbb{Z}: \(\frac{x}{p}\)=-1\bigg\}.$$
There are explicit formulas for the sums
$$
A_p:=\sum_{0<x<\frac{p}{2},x\in\mathcal{R}}x.
$$
If we put
$$
B_p:=\sum_{0<x<\frac{p}{2},x\in\mathcal{N}}x,
$$
then in view of
\begin{equation}\label{ap+bp}
A_p+B_p=1+2+\cdots+\frac{p-1}{2}=\frac{p^2-1}{8}
\end{equation}
it suffices to evaluate $A_p-B_p$. In case $p\equiv3\pmod4$  we have
\begin{equation}\label{ap-bp}
A_p-B_p=\begin{cases}0&\mbox{if}\ p\equiv7\pmod8,
\\ph(-p)&\mbox{if}\ p\equiv3\pmod8,\end{cases}
\end{equation}
which is a special case of the equation (23) in Lerch's paper \cite{Lerch}.  To state the result for
the case $p\equiv1\pmod4$ we need to introduce  Dirichlet L-series.
Let
$$L(s,(\cdot/p))=\sum_{n>0}\(\frac{n}{p}\)n^{-s}$$
be Dirichlet L-series attached to the Legendre symbol $(\frac{\cdot}{p})$, where $s$ is any
complex number with {\rm Re} $s>1$. It is well known that
$L(s,(\cdot/p))$  can be analytically continued to the whole complex
plane. In the remaining part of this paper, we assume that $L(s,(\cdot/p))$ is defined on the whole
complex plane. By Berndt's survey paper \cite{BE76} we know that if we let $\chi$ denote the Legendre
symbol for $p$, then
\begin{equation}\label{2ap-bp}
A_p-B_p= -\frac{\tau(\chi)p}{\pi^2}  (1-\chi(2)/4) L(2,\chi),
\end{equation}
where
$$
L(2,\chi)=-\frac{\tau(\chi)\pi^2}{p^2}B_{2,\chi}
$$
and
$$
B_{2,\chi}=-2L(-1,\chi).
$$
By (\ref{ap+bp})--(\ref{2ap-bp}) one can obtain the following results:
\begin{itemize}

\item If $p\equiv3\pmod4$, then
\begin{equation}\label{ap3}
A_p=\begin{cases}(p^2-1)/16&\mbox{if}\ p\equiv7\pmod8,
\\(p^2-1+8ph(-p))/16&\mbox{if}\ p\equiv3\pmod8,\end{cases}
\end{equation}
where $h(-p)$ is the class number of $\mathbb{Q}(\sqrt{-p})$.

\item  It $p\equiv1\pmod4$, then
\begin{equation}\label{ap4}
A_p=\begin{cases}(p^2-1+12\cdot L((\cdot/p),-1))/16&\mbox{if}\ p\equiv1\pmod8,
\\(p^2-1+20\cdot L((\cdot/p),-1))/16&\mbox{if}\ p\equiv5\pmod8.\end{cases}
\end{equation}
\end{itemize}

Motivated by the above rich results on quadratic residues, in this paper we concentrate on the
products concerning quadratic residues.  Let $p$ be an odd prime and set
 $$M_p:=\prod_{0<x<p/2,x\in\mathcal{R}}x.$$
We will determine $M_p$\ {\rm mod}\ $p$.
It turns out that these results involve quartic Gauss sums, Stickelberger's congruence and values of
Dirichlet L-series at negative integers. By Wilson's theorem one may easily verify that
$$M_p^2\equiv\begin{cases}1\pmod p&\mbox{if}\ p\equiv 5\pmod8,\\-1\pmod p&\mbox{if}\ p\equiv
1\pmod8,\end{cases}$$
and hence we have
$$M_p\equiv\begin{cases}\pm1\pmod p&\mbox{if}\ p\equiv 5\pmod8,\\\pm\frac{p-1}{2}!\pmod p&\mbox{if}\
p\equiv
1\pmod8.\end{cases}$$
We shall show that the value of  $M_p$\ {\rm mod} $p$  behaves quite differently according to
$p\equiv5\pmod8$ or $p\equiv1\pmod8$. So we discuss the two cases separately.

We first consider the case $p\equiv5\pmod8$.  To state our result we need to introduce the rational
$4$-th power residue symbol.  For any prime $p\equiv1\pmod4$  and  integer $a$ we define the rational
$4$-th power residue symbol as follows:
$$\(\frac{a}{p}\)_4=\begin{cases}0&\mbox{if}\ p\mid a,
\\1&\mbox{if}\ \text{{\it a} is a 4-th power residue modulo {\it p}} ,
\\-1&\mbox{otherwise}.\end{cases}$$
Now we state our result for the case $p\equiv5\pmod8$.
\begin{theorem}\label{Theorem products of quadratic residues p=5 mod8}
Let $p\equiv5\pmod8$ be a prime. Then we have
$$\prod_{0<x<p/2,x\in\mathcal{R}}x\equiv(-1)^{1+\#\{0<x<p/2:\ (\frac{x}{p})_4=1\}}\pmod p.$$
\end{theorem}
We now consider the case $p\equiv1\pmod8$. This case is related to quartic Gauss sums and
Stickelberger's congruence. Here we give a brief review of these (for more details readers may refer
to \cite[Chapter 3]{Cohen1}). Let $p\equiv1\pmod8$ be a prime, and let $\zeta_{p-1}=e^{2\pi {\bf
i}/(p-1)}$. Let $K=\mathbb{Q}(\zeta_{p-1},\zeta_p)$, and let $\mathcal{O}_K$ be the ring of algebraic
integers of $K$. Let $\mathfrak{p}$ be a prime ideal of $\mathcal{O}_K$ lying above the prime ideal
$(1-\zeta_p)\mathbb{Z}[\zeta_p]$ of $\mathbb{Z}[\zeta_p]$.
Clearly we have $$\mathcal{O}_K/\mathfrak{p}\cong\mathbb{Z}/p\mathbb{Z}=\mathbb{F}_p, $$
where $\mathbb{F}_p$ denotes the finite field with $p$ elements. It is known that the map
 $$u_{\mathfrak{p}}: \zeta_p^k\mapsto\zeta_p^k\pmod{\mathfrak{p}}$$
  is a bijection from $\{\zeta_{p-1}^k: k=0,1,\cdots,p-2\}$ onto
  $(\mathcal{O}_K/\mathfrak{p})^{\times}\cong(\mathbb{Z}/p\mathbb{Z})^{\times}=\mathbb{F}_p^{\times}$,
where $R^{\times}$ denotes the group of invertible elements of the ring $R$. We now let
$\omega_{\mathfrak{p}}=u_{\mathfrak{p}}^{-1}$ be the mulplicative character of $\mathbb{F}_p$. Clearly
$\omega_{\mathfrak{p}}$ generates the character group $\chi(\mathbb{F}_p)$. We let
$\chi_{\pi}=\omega_{\mathfrak{p}}^{-(p-1)/4}$ be the character of order $4$. We consider the Jacobi
sum
$J(\chi_{\pi},\chi_{\pi})$. By \cite[Proposition 3.6.4]{Cohen1} we have
\begin{equation}\label{Eq. Congruence of Jacobi sums}
J(\chi_{\pi},\chi_{\pi})\equiv -\frac{\frac{p-1}{2}!}{(\frac{p-1}{4}!)^2}\pmod{\mathfrak{p}}.
\end{equation}
Moreover, it is known that (cf. \cite[Theorem 3.9]{BE79})
\begin{equation}\label{Eq. value of Jacobi sums}
J(\chi_{\pi},\chi_{\pi})=a+4b{\bf i}
\end{equation}
with $a\equiv-1\pmod4$ and $p=a^2+16b^2$.
Now we consider the Gauss sum $G(\chi_{\pi})$. In 1977, Loxton \cite{Loxton} posed the following
conjecture concerning the explicit value of $G(\chi_{\pi})$ (in a slightly different form):
\begin{equation}\label{Eq. Gauss sums}
G(\chi_{\pi})=C_p\(\frac{\vert b\vert}{ \vert a\vert}\)(-1)^bp^{1/4}J(\chi_{\pi},\chi_{\pi})^{1/2},
\end{equation}
where $(\frac{\cdot}{\vert a\vert})$ is the Jacobi symbol, $p^{1/4},\ {\rm Re}\
J(\chi_{\pi},\chi_{\pi})^{1/2}>0$ and $C_p$ is defined by
\begin{equation}\label{Eq. C}
C_p=\pm1\ \text{and}\ C_p\equiv \frac{4\vert b\vert}{a}\(\frac{p-1}{2}\)!\pmod p.
\end{equation}
Later Matthews \cite{Matthews} confirmed Loxton's elegant conjecture. Hence almost 175
years after Gauss's determination of quadratic Gauss sums, there has been found an elegant
formula for quartic Gauss sum. Readers may refer to \cite{BE} for the history and details on Gauss
sums.

Recall that $\ve_p=\frac{u_p+v_p\sqrt{p}}{2}>1$ ($u_p,v_p\in\mathbb{Z}$) and $h(p)$ are the
fundamental unit and the class number of $\mathbb{Q}(\sqrt{p})$ respectively. Now we consider the
number field $L=K(\ve_p^{1/2},J(\chi_{\pi},\chi_{\pi})^{1/2})$. Let $\mathcal{O}_L$ be the ring of
algebraic integers of $L$. Let $\mathfrak{P}$ be a prime ideal of $\mathcal{O}_L$ lying above
$\mathfrak{p}$. By \cite[Corollary 1.1]{Sun ffa1} if we write $\ve_p^{h(p)}=a_p+b_p\sqrt{p}$ with
$2a_p,2b_p$ integers, then we have
\begin{equation}\label{Eq. Congruence of fundamental unit}
a_p\equiv-\frac{p-1}{2}!\pmod p.
\end{equation}
For more results on the congruences involving fundamental units, readers may refer to
\cite{Chowla,Wu-She}. Combining Eq. (\ref{Eq. Congruence of Jacobi sums}) with Eq. (\ref{Eq.
Congruence of fundamental unit}), we obtain
$$\frac{\ve_p^{h(p)}}{(((p-1)/4)!)^2\cdot J(\chi_{\pi},\chi_{\pi})}\equiv 1\pmod {\mathfrak{p}}.$$
We therefore define $\beta_p$ by
\begin{equation}\label{Eq. beta p}
\beta_p=0,1,\ \text{and}\ (-1)^{\beta_p}\equiv \frac{\ve_p^{h(p)/2}}{((p-1)/4)!\cdot
J(\chi_{\pi},\chi_{\pi})^{1/2}}\pmod{\mathfrak{P}},
\end{equation}
where $\ve_p^{1/2}$ and Re $J(\chi_{\pi},\chi_{\pi})^{1/2}>0$.
Now we state our  result for the case $p\equiv1\pmod8$ .
\begin{theorem}\label{Theorem products of quadratic residues p=1 mod8}
Let $p\equiv1\pmod8$ be a prime. If we write $p=a^2+16b^2$ with $a,b\in\mathbb{Z}$. Then we have
$$\prod_{0<x<p/2,x\in\mathcal{R}}x\equiv C_p(-1)^{1+\beta_p+\lfloor\frac{p}{8}\rfloor}\(\frac{\vert
b\vert}{\vert a\vert}\)\(\frac{p-1}{2}!\)\pmod p,$$
where $C_p=\pm1$ is defined by {\rm Eq. (\ref{Eq. C})} and $(\frac{\cdot}{\vert a\vert})$ is the
Jacobi symbol.
\end{theorem}
We will prove our main results in Section 2.  In Section 3 we  pose some problems for further
research.
\maketitle
\section{Proofs of the main results}
In this section, we also adopt the notations in Section 1. We first prove Theorem \ref{Theorem
products of quadratic residues p=5 mod8}.

{\bf Proof of Theorem \ref{Theorem products of quadratic residues p=5 mod8}.}
Let $p\equiv5\pmod8$ be a prime. Clearly we have
$$-1\cdot\prod_{0<x<p/2,x\in\mathcal{R}}x^2\equiv
\prod_{0<x<p/2,x\in\mathcal{R}}x(p-x)\equiv\prod_{0<x<p,x\in\mathcal{R}}x\equiv-1\pmod p.$$
Hence $\prod_{0<x<p/2,x\in\mathcal{R}}x\equiv\pm1\pmod p$. If we set
$$r=\#\bigg\{0<x<p/2: \(\frac{x}{p}\)_4=1\bigg\},$$
then we have
$$\prod_{0<x<p/2,x\in\mathcal{R}}\(\frac{x}{p}\)_4
=\(\frac{\prod_{0<x<p/2,x\in\mathcal{R}}x}{p}\)_4=(-1)^{\frac{p-1}{4}-r}.$$
Noting that $-1$ is a $4$-th non-residue modulo $p$, we therefore have
$$\prod_{0<x<p/2,x\in\mathcal{R}}x\equiv (-1)^{\frac{p-1}{4}-r}=(-1)^{r+1}\pmod p.$$
This completes the proof.\qed

To prove Theorem \ref{Theorem products of quadratic residues p=1 mod8} we need to know the  explicit
value of the following product:
$$\prod_{0<x<p/2,x\in\mathcal{R}}(1-\zeta_p^{2x}).$$
Let $\zeta_p=e^{2\pi{\bf i}/p}$. It is clear that for any integer $k$ with $p\nmid k$ we have
 $$(1-\zeta_p^{k})/(1-\zeta_p)  \equiv k\pmod{(1-\zeta_p)\mathbb{Z}[\zeta_p]}.$$
 One deduces that
$$2^{\frac{p-1}{4}}\times\prod_{0<x<p/2,x\in\mathcal{R}}x\equiv
\prod_{0<x<p/2,x\in\mathcal{R}}\frac{1-\zeta_p^{2x}}{1-\zeta_p}\
\pmod{(1-\zeta_p)\mathbb{Z}[\zeta_p]}.$$
Defining
$$W_p:=\prod_{0<x<p/2,x\in\mathcal{R}}(1-\zeta_p^{2x})$$
gives that
$$2^{\frac{p-1}{4}}\times\prod_{0<x<p/2,x\in\mathcal{R}}x\equiv \frac{W_p}{(1-\zeta_p)^{(p-1)/4}}\
\pmod{(1-\zeta_p)\mathbb{Z}[\zeta_p]}.$$
The following lemma gives the explicit value of $W_p$.
\begin{lemma}\label{Theorem A}
Let $p\equiv1\pmod4$ be a prime. Then we have
$$\prod_{0<x<\frac{p}{2},x\in\mathcal{R}}(1-\zeta_p^{2x})=\begin{cases}
(-1)^{\lfloor\frac{p}{8}\rfloor}\zeta_p^{\frac{p^2-1+12L(-1,(\cdot/p))}{16}}p^{1/4}\ve_p^{-h(p)/2}&\mbox{if}\
8\mid p-1,
\\\\(-1)^{1+\lfloor\frac{p}{8}\rfloor}\zeta_p^{\frac{p^2-1+20L(-1,(\cdot/p))}{16}}{\bf i}\cdot
p^{1/4}\ve_p^{h(p)/2}&\mbox{if}\ 8\mid p-5,\end{cases}$$
\end{lemma}
where $\lfloor\cdot\rfloor$ is the floor function and $p^{1/4}, \ve_p^{1/2}>0$.
\begin{proof}
Let's still set
$$W_p:=\prod_{0<x<\frac{p}{2},x\in\mathcal{R}}(1-\zeta_p^{2x}).$$
We first compute the absolute value of $W_p$. By definition we have
\begin{equation}\label{Eq. Absolute value of Wp}
W_p\cdot
\overline{W_p}=\prod_{k=1}^{(p-1)/2}(1-\zeta_p^{2k^2})=\sqrt{p}\cdot\ve_p^{-(\frac{2}{p})h(p)}.
\end{equation}
The last equality is a known result (cf. \cite[Theorem 1.3]{Sun ffa}) and $\overline{W_p}$ denotes the
conjugation of $W_p$. Now we determine the argument
$\Arg(W_p)$ of the complex number $W_p$. We have the following equalities:
\begin{align}\label{Eq. Wp}
W_p=\prod_{0<x<\frac{p}{2},x\in\mathcal{R}}-\zeta_p^x(\zeta_p^x-\zeta_p^{-x})
=(-1)^{\frac{p-1}{4}}\cdot\zeta_p^{A_p}\prod_{0<x<\frac{p}{2},x\in\mathcal{R}}2{\bf i}\sin \frac{2\pi
x}{p},
\end{align}
where $A_p$ is defined as in the introduction. From this we obtain
$$\Arg(W_p)\equiv\begin{cases}\lfloor p/8\rfloor\pi+2\pi A_p/p\pmod{2\pi\mathbb{Z}}&\mbox{if}\
p\equiv1\pmod8,
\\-\pi/2+\lfloor p/8\rfloor\pi+2\pi A_p/p\pmod{2\pi\mathbb{Z}}&\mbox{if}\ p\equiv5\pmod8,\end{cases}$$
where $\lfloor\cdot\rfloor$ is the floor function.  By the explicit formula for $A_p$ in the
introduction we have
$$W_p=\begin{cases}
(-1)^{\lfloor\frac{p}{8}\rfloor}\zeta_p^{\frac{p^2-1+12L(-1,(\cdot/p))}{16}}p^{1/4}\ve_p^{-h(p)/2}&\mbox{if}\
p\equiv1\pmod8,
\\\\(-1)^{1+\lfloor\frac{p}{8}\rfloor}\zeta_p^{\frac{p^2-1+20L(-1,(\cdot/p))}{16}}{\bf i}\cdot
p^{1/4}\ve_p^{h(p)/2}&\mbox{if}\ p\equiv5\pmod8.\end{cases}$$
This completes the proof.
\end{proof}
To prove Theorem \ref{Theorem products of quadratic residues p=1 mod8} we also need  the following
lemma,  which is known as Stickelberger's congruence (cf. \cite[Theorem 3.6.6]{Cohen1}).
\begin{lemma}\label{Lemma s congruence}
Let notations be as in Section 1, and let $p\equiv1\pmod8$ be a prime. Then for all $0\le r\le p-2$ we
have
$$\frac{G(\omega_{\mathfrak{p}}^{-r})}{(\zeta_p-1)^r}\equiv \frac{-1}{r!}\pmod {\mathfrak{p}},$$
where $G(\omega_{\mathfrak{p}}^{-r})$ is the Gauss sum associated to the character
$\omega_{\mathfrak{p}}^{-r}$.
\end{lemma}
Now we are in a position to prove our last theorem.

{\bf Proof of Theorem \ref{Theorem products of quadratic residues p=1 mod8}.}
By Lemma \ref{Lemma s congruence} we have
$$\frac{G(\chi_{\pi})}{(1-\zeta_p)^{(p-1)/4}}\equiv\frac{-1}{((p-1)/4)!}\pmod {\mathfrak{P}}.$$
By Eq. (\ref{Eq. Gauss sums}) we obtain
$$\frac{C_p(\frac{\vert b\vert}{\vert
a\vert})(-1)^bp^{1/4}J(\chi_{\pi},\chi_{\pi})^{1/2}}{(1-\zeta_p)^{(p-1)/4}}
\equiv\frac{-1}{((p-1)/4)!}\pmod{\mathfrak{P}}.$$
If we write $p=a^2+16b^2$ with $a,b\in\mathbb{Z}$, then it is known that $2$ is a $4$-th residue
modulo $p$ if and only if $2\mid b$. Hence we have
$$\frac{p^{1/4}}{(1-\zeta_p)^{(p-1)/4}}\equiv\frac{-C_p(\frac{\vert b\vert}{\vert
a\vert})(\frac{2}{p})_4}{((p-1)/4)!\cdot J(\chi_{\pi},\chi_{\pi})^{1/2}}\pmod {\mathfrak{P}}.$$
Combining this with Lemma \ref{Theorem A}, we obtain
\begin{align*}
\prod_{0<x<p/2,x\in\mathcal{R}}x\equiv
\frac{(-1)^{1+\lfloor\frac{p}{8}\rfloor}\ve_p^{-h(p)/2}C_p(\frac{\vert b\vert}{\vert
a\vert})}{((p-1)/4)!\cdot J(\chi_{\pi},\chi_{\pi})^{1/2}}\pmod {\mathfrak{P}}.
\end{align*}
By Eq. (\ref{Eq. Congruence of fundamental unit}), we have $\ve_p^{h(p)}\equiv-\frac{p-1}{2}!\pmod
{\mathfrak{P}}$. From this, we finally obtain that
\begin{align*}
\prod_{0<x<p/2,x\in\mathcal{R}}x&\equiv
\frac{(-1)^{1+\lfloor\frac{p}{8}\rfloor}\ve_p^{h(p)/2}C_p(\frac{\vert b\vert}{\vert
a\vert})}{((p-1)/4)!\cdot J(\chi_{\pi},\chi_{\pi})^{1/2}}\cdot\(\frac{p-1}{2}!\)\\
&\equiv C_p(-1)^{1+\beta_p+\lfloor\frac{p}{8}\rfloor}\(\frac{\vert b\vert}{\vert
a\vert}\)\(\frac{p-1}{2}!\)\pmod {\mathfrak{P}}.
\end{align*}
This implies our desired result.\qed
\maketitle
\section{Some open problems}
In this section, we pose some open problems for further research. For any $k,n\in \mathbb{Z}^+$, we
define
$$H_k^{(n)}:=\sum_{1\le x\le k}\frac{1}{x^n}.$$
These numbers are known as harmonic numbers. Let $p>3$ be a prime. In 1862 Wolstenholme showed that
$$H_{p-1}^{(1)}\equiv0\pmod{p^2}.$$
Later Lehmer \cite{L} determined $H_{\frac{p-1}{2}}^{(1)}$ {\rm mod} $p^2$ completely. Motivated by
the above work, we define
$$H_{\mathcal{R},\frac{p-1}{2}}^{(n)}:=\sum_{0<x<p/2,x\in\mathcal{R}}\frac{1}{x^n}.$$
It is easy to see that if $p\equiv3\pmod4$, then
$$H_{\mathcal{R},\frac{p-1}{2}}^{(1)}\equiv \frac{1}{2}H_{\frac{p-1}{2}}^{(1)}\pmod p,$$
and that if $p\equiv1\pmod4$, then
$$H_{\mathcal{R},\frac{p-1}{2}}^{(2)}\equiv0\pmod p.$$

In view of the above, we now pose our first problem:
\begin{problem}
  {\rm (i)  }Let $p\equiv1\pmod4$ be a prime.  Can we determine the explicit value of
  $H_{\mathcal{R},\frac{p-1}{2}}^{(1)}$ {\rm mod} $p$ ?

{\rm (ii)  } Let $p\equiv3\pmod4$ be a prime. Can we determine the explicit value of
$H_{\mathcal{R},\frac{p-1}{2}}^{(2)}$ {\rm mod} $p$ ?
\end{problem}

The following table gives the  values of $H_{\mathcal{R},\frac{p-1}{2}}^{(1)}$ {\rm mod} $p$ for
primes less than $100$ with $p\equiv1\pmod4$:
\begin{center}
\begin{tabular}{|c|c|c|c|c|c|c|c|c|c|c|c|}

\hline
 p&5&13&17&29&37&41&53&61&73&89&97\\
\hline
$H_{\mathcal{R},\frac{p-1}{2}}^{(1)}$ {\rm mod} $p$ &1&7&4&23&12&18&10&13&17&83&40\\
\hline
\end{tabular}
\end{center}

We also calculated  the  values of $H_{\mathcal{R},\frac{p-1}{2}}^{(2)}$ {\rm mod} $p$ for  primes
less than $100$ with $p\equiv3\pmod4$:
\begin{center}
\begin{tabular}{|c|c|c|c|c|c|c|c|c|c|c|c|c|c|}

\hline
 p&3&7&11&19&23&31&43&47&59&67&71&79&83\\
\hline
$H_{\mathcal{R},\frac{p-1}{2}}^{(2)}$ {\rm mod} $p$ &1&3&8&5&19&13&29&17&14&18&56&40&14\\
\hline
\end{tabular}
\end{center}

Now we turn to a problem concerning the products of quadratic residues.
\begin{problem}
   Let $p\equiv3\pmod4$ be a prime. Let $M_p=\prod_{0<x<p/2,x\in\mathcal{R}}x$.  Is there any pattern
   for the values of $M_p$ {\rm mod} $p$ ?
\end{problem}
The following table gives the  values of $M_p$ {\rm mod} $p$ for  primes less than $100$ with
$p\equiv3\pmod4$:
\begin{center}
\begin{tabular}{|c|c|c|c|c|c|c|c|c|c|c|c|c|c|}

\hline
 p&3&7&11&19&23&31&43&47&59&67&71&79&83\\
\hline
$M_p$ {\rm mod} $p$ &1&2&5&17&18&5&41&4&29&10&58&38&51\\
\hline
\end{tabular}
\end{center}

\subsection*{Acknowledgements}
We are grateful to the referee for a careful reading of the original manuscript and for helpful
comments
which improved the quality of our paper. This research is supported by the National Natural Science
Foundation of China (Grant No. 11971222). The first author is also supported by NUPTSF (Grant No.
NY220159).

\end{document}